\documentclass[11 pt,reqno]{amsart}

\usepackage[latin1]{inputenc}
\usepackage{amsmath,amsthm,amssymb,amscd,mathrsfs,mathtools,xcolor,esint}
\usepackage[margin=1in]{geometry}
\usepackage[shortlabels]{enumitem}
\definecolor{darkgreen}{rgb}{0.0, 0.6, 0.13}
\usepackage{hyperref} 
\usepackage{macros_for_trees}

\usepackage{scalerel}

\usepackage{tikz}
\usetikzlibrary{fit,positioning}
\usepackage{pgfplots}
\pgfplotsset{width=7cm,compat=1.8}
\usepackage{pgfplotstable}

\usepackage{soul}
\definecolor{airforceblue}{rgb}{0.0, 0.30, 0.69}

\newtheorem{thm}{Theorem}[section]
 
 \newtheorem{lem}[thm]{Lemma}

 \theoremstyle{definition}
 
 \theoremstyle{remark}
 \newtheorem{rem}[thm]{Remark}
 
 \numberwithin{equation}{section}

\usepackage{multirow}
\usepackage{array}
\usepackage{caption}
\usepackage{subcaption}
\usepackage{booktabs} 
\usepackage{tablefootnote} 
\usepackage{xcolor,colortbl}

 \makeatletter
\def\@tocline#1#2#3#4#5#6#7{\relax
  \ifnum #1>\c@tocdepth 
  \else
    \par \addpenalty\@secpenalty\addvspace{#2}
    \begingroup \hyphenpenalty\@M
    \@ifempty{#4}{
      \@tempdima\csname r@tocindent\number#1\endcsname\relax
    }{
      \@tempdima#4\relax
    }
    \parindent\z@ \leftskip#3\relax \advance\leftskip\@tempdima\relax
    \rightskip\@pnumwidth plus4em \parfillskip-\@pnumwidth
    #5\leavevmode\hskip-\@tempdima
      \ifcase #1
       \or\or \hskip 1em \or \hskip 2em \else \hskip 3em \fi
      #6\nobreak\relax
    \hfill\hbox to\@pnumwidth{\@tocpagenum{#7}}\par
    \nobreak
    \endgroup
  \fi}
\makeatother

\begin{document}
\author[Yu Deng]{Yu Deng$^1$}
\address{$^1$ Department of Mathematics, University of Southern California, Los Angeles,  CA 90089, USA }
\email{yudeng@usc.edu}
\thanks{$^1$Y. D. is funded in part by NSF-DMS-2246908 and a Sloan Fellowship.}
\author[Andrea R. Nahmod]{Andrea R. Nahmod$^2$}
\address{$^2$ 
Department of Mathematics,  University of Massachusetts,  Amherst MA 01003}
\email{nahmod@math.umass.edu}
\thanks{$^2$A.N. is funded in part by NSF-DMS-2052740, NSF-DMS-2101381 and the Simons Foundation Collaboration Grant on Wave Turbulence (Nahmod's Award ID 651469).}
\author[Haitian Yue]{Haitian Yue$^3$}
\address{$^3$ Institute of Mathematical Sciences, ShanghaiTech University, Shanghai, 201210, China}
\email{yuehaitian@shanghaitech.edu.cn}
\thanks{$^3$H.Y. is funded in part by the Shanghai Technology Innovation
Action Plan (No.22JC1402400) and a Chinese overseas high-level young talents program (2022).}
  \date{}
 \title{The Probabilistic Scaling Paradigm}
 \dedicatory{Dedicated to Carlos E. Kenig}

 \begin{abstract} In this note we further discuss the probabilistic scaling introduced by the authors in \cite{DNY19,DNY20}. In particular we do a case study comparing the stochastic heat equation, the nonlinear wave equation and the nonlinear Schr\"odinger equation.

 \bigskip
 \bigskip

\noindent \textbf{Keywords}: Random data theory, probabilistic scaling, NLS, NLW, stochastic heat equation, probabilistic well posedness; Gibbs measure.
\bigskip

\noindent \textbf{Mathematics Subject Classification (2020)}:  35Q60, 35R60 (15A69 15B52 37E20 60H30).
 
 \end{abstract}
 \maketitle
 
\tableofcontents
 \section{Introduction}\label{intro} 
 
In recent years we have been interested in the  general question of   understanding and describing how randomness affects the behavior of solutions to PDEs, particularly the study of propagation of randomness.
Randomness may come into the problem in various ways but two common ones are: 
from the equation such as in stochastic problems with additive or multiplicative noise;
 from random initial data which obeys some canonical law of distribution (e.g. Gaussian law).

Some{\footnote{{{One may also ask about how their distribution evolve over time; at least how do the ensemble averages evolve. Such questions arise in wave turbulence; see \cite{BGHS21}; \cite{DH, DH2, DH2', DH3}; \cite{CG, CG20} and references therein.}}}} fundamental questions are:
  {what is the optimal regime where the solution exists and is unique almost surely, at least locally in time;}
 {can one describe the solution in terms of the random structure of the initial data at least for short times;}
 {if  there are (formally invariant) Gibbs measures, can we justify their invariance.}
In these works \cite{DNY19, DNY20,DNY21} and in subsequent joint work with B. Bringmann  \cite{BDNY22} we provided some answers to these questions in the context of nonlinear Schr\"odinger equations (NLS) and nonlinear wave equations (NLW).

In this paper we reprise and extend a cornerstone idea that lies behind the works above. We
consider general evolution PDEs of parabolic, hyperbolic and dispersive type, in the probabilistic setting. We are interested in the local theory of these equations in different regularity spaces, and in particular the threshold between probabilistic almost-sure well-posedness and ill-posedness.

 We will discuss the \emph{scaling heuristics} for each type of equations that predicts a \emph{critical} regularity exponent. This is a guiding principle for the study of probabilistic local well-posedness, and in many (but not all cases this indeed matches the actual local well-posedness and ill-posedness results. In fact, in the parabolic and deterministic settings, it is not uncommon that the scaling does not match optimal local well-posedness threshold\footnote{Heuristically, one can check the following example  \cite{Weber}: $\dot{x}= \sigma(x)\ \dot{w}$. At least when $w$ is fractional Brownian motion, there is no canonical stochastic definition of iterated integrals below the threshold $\gamma = \frac{1}{4}$, where $\gamma$ refers the regularity of $w$. However, the regime $\gamma>0$ is in principle  subcritical in the parabolic setting.}, but the scaling criticality still acts as a good guiding principle.

\smallskip
In the \emph{parabolic} setting, the notion of scaling has been systematically studied in earlier works \cite{Hairer, Hairer3}. By using the theory of regularity structures \cite{Hairer}, it has been proved that, under suitable additional restrictions and with suitable renormalization, any nonlinear parabolic SPDE that is \emph{above} the scaling threshold must be almost-surely locally well-posed.

In the \emph{hyperbolic} and \emph{dispersive} settings, the situation is much more complicated. In \cite{DNY19,DNY20}, the authors introduced the notion of \emph{probabilistic scaling} for nonlinear Schr\"{o}dinger equations, and proved the analog result to \cite{Hairer}, namely that any nonlinear Schr\"{o}dinger equation (of at least cubic nonlinearity) that is \emph{above} the scaling threshold must be almost-surely locally well-posed.

However, unlike parabolic equations, the behavior of dispersive equations, especially that of \emph{resonances}, is much richer and contains more possibilities. For this reason, results similar to \cite{DNY20} have not been obtained for general dispersion relations.

\smallskip
In this note, we will describe a general framework of the heuristic (probabilistic) scaling argument that treats the cases of parabolic, hyperbolic and dispersive equations in a unified manner. In addition, we will describe the cases in which this scaling heuristics is known or expected to be precise, as well as possible reasons that may cause it to deviate from the actual threshold. Finally, we conclude by connecting this scaling heuristics to recent advancements in wave turbulence theory.

\section{The scaling heuristics}\label{scaling} In this section we describe the unified scaling argument that covers parabolic, hyperbolic and dispersive equations. For the purpose of demonstration, we will fix the domain of our equation to be the $d$ dimensional torus $\mathbb{T}^d$, and consider the following three model equations:
\begin{itemize}
\item The stochastic heat equation:
\begin{equation}\label{she}(\partial_t-\Delta)u+\mathcal{N}_p(u)=\xi,
\end{equation} where $\xi$ is some spacetime Gaussian noise on $\mathbb{R}\times\mathbb{T}^d$;
\item The semilinear wave equation:
\begin{equation}\label{nlw}(\partial_t^2-\Delta)u+\mathcal{N}_p(u)=0,
\end{equation} with initial data being a Gaussian random field on $\mathbb{T}^d$.
\item The semilinear Schr\"{o}dinger equation:
\begin{equation}\label{nls}(i\partial_t-\Delta)u+\mathcal{N}_p(u)=0,
\end{equation} with initial data being a Gaussian random field on $\mathbb{T}^d$.
\end{itemize}
Here in (\ref{she})--(\ref{nls}), $\mathcal{N}_p(u)$ is a $p$-th order polynomial of $u$ without derivative (the exact form does not matter in local theory), and the precise form of noise or initial data will be specified in Section \ref{data} below.
\subsection{Choice of noise and initial data}\label{data}
\subsubsection{The parabolic case} We start with the heat equation (\ref{she}). A canonical choice of the noise would be a fractional derivative (or antiderivative) of spacetime white noise, namely
\begin{equation}\label{data1} \xi(t,x)=\sum_{k\in\mathbb{Z}^d}\langle k\rangle^{-\alpha}e^{ik\cdot x}\xi_k(t),
\end{equation} for some $\alpha\in\mathbb{R}$, where $\xi_k=\mathrm{d}B_k/\mathrm{d}t$ are derivatives of independent one-dimensional Brownian motions. In terms of regularity, we know that
\[\xi\in C_t^0C_x^{\alpha-1-(d/2)-\varepsilon}\] almost surely\footnote{Here and also in the following parts of this paper we write $C^\gamma$ for the Besov space $B_{\infty,\infty}^{\gamma}$.}. Note that the linear noise component of the solution $u$ to (\ref{she}) is given by
\begin{equation}\label{linearnoise}u_{\mathrm{lin}}:=\int_0^t e^{(s-t)\Delta}\eta(s)\,\mathrm{d}s;
\end{equation} since the parabolic regularity gains two derivatives, it is easy to see that
\begin{equation}\label{regdata1} u_{\mathrm{lin}}\in C_t^0C_x^{\alpha+1-(d/2)-\varepsilon}
\end{equation} almost surely, which is the space in which the solution $u$ is constructed provided the initial data $u(0)$ has matching regularity (see for example \cite{Hairer}). Ignoring the arbitrary small constant $\varepsilon$, we may call the value $\alpha+1-(d/2)$ the \emph{regularity level} of the solution  to (\ref{she}) with noise (\ref{data1}).
\subsubsection{The hyperbolic case} We now turn to the wave equation (\ref{nlw}), for which the randomness comes from the initial data. We fix the initial data $(u_0,u_1)=(u(0),\partial_tu(0))$ as
\begin{equation}\label{data2}u_0=\sum_{k\in\mathbb{Z}^d}\langle k\rangle^{-\alpha-1}e^{ik\cdot x}g_k,\quad  u_1=\sum_{k\in\mathbb{Z}^d}\langle k\rangle^{-\alpha}e^{ik\cdot x}g_k',
\end{equation} where $g_k$ and $g_k'$ are i.i.d. Gaussians, for some $\alpha\in\mathbb{R}$. In terms of regularity, we have that
\[u_0\in C_x^{\alpha+1-(d/2)-\varepsilon},\quad u_1\in C_x^{\alpha-(d/2)-\varepsilon}\] almost surely, so the linear component of the solution
\begin{equation}\label{lineardata}u_{\mathrm{lin}}:=\cos(t|\nabla|)u_0+\frac{\sin(t|\nabla|)}{|\nabla|}u_1\end{equation} satisfies that
\begin{equation}\label{datagre2}
u_{\mathrm{lin}}\in C_x^{\alpha+1-(d/2)-\varepsilon}
\end{equation} almost surely. In fact, by the properties of the i.i.d. Gaussians, it is easy to see that  $u_{\mathrm{lin}}$ is also in $H_x^{\alpha+1-(d/2)-\varepsilon}$ almost surely. Similar to above, we call the value $\alpha+1-(d/2)$ the \emph{regularity level} of the solution  to (\ref{nlw}) with data (\ref{data2}).
\subsubsection{The dispersive case} The case of the Schr\"{o}dinger equation (\ref{nls}) is similar to (\ref{nlw}), but we only need to specify the initial data $u(0)$. Suppose
\begin{equation}\label{data3}
u(0)=u_0=\sum_{k\in\mathbb{Z}^d}\langle k\rangle^{-\alpha-1}e^{ik\cdot x}g_k,
\end{equation} where $g_k$ are i.i.d. Gaussians, then in the same way as (\ref{nlw}), we have
\[u_0\in C_x^{\alpha+1-(d/2)-\varepsilon},\] and hence the linear component of the solution \begin{equation}\label{lineardata}u_{\mathrm{lin}}:=e^{-it\Delta}u_0\end{equation} satisfies that
\begin{equation}\label{datagre2}
u_{\mathrm{lin}}\in C_x^{\alpha+1-(d/2)-\varepsilon}
\end{equation} almost surely. Similarly as in the hyperbolic case  $u_{\mathrm{lin}}$ is also in $H_x^{\alpha+1-(d/2)-\varepsilon}$ almost surely.  Again we call the value $\alpha+1-(d/2)$ the \emph{regularity level} of the solution  to (\ref{nls}) with data (\ref{data3}).
\subsection{The criticality threshold}\label{sec:criticality} Note that the noise (\ref{data1}), as well as the initial data (\ref{data2}) and (\ref{data3}), are essentially  {homogeneous} when restricted to high frequency/small physical scales; in other words, rescaling the noise (\ref{data1}) or the initial data (\ref{data2}) or (\ref{data3}) just corresponds to restricting to frequency $\sim N$ for some $N\gg 1$.

In general, we can view the criticality of a given equation in a given space $C^s$ (or $H^s$) in the following aspect: suppose one starts with some initial data $f$, which is rescaled to frequency $\sim N\gg 1$ and normalized to have $C^s$ norm $\sim 1$, and computes the second order iterations of the nonlinearity with this $f$. Then, the equation is \emph{subrcritical} if the result of iteration carries a \emph{negative} power of $N$ when measured in $C^s$ norm, is \emph{critical} if the result of iteration is \emph{independent} of $N$ when measured in $C^s$ norm, and is \emph{supercritical} if the result of iteration carries a \emph{positive} power of $N$ when measured in $C^s$ norm. Due to the above scaling covariant observation, it is easy to see that the above definition is equivalent to the usual definition of scaling criticality using rescaling.

The advantage of such definition of criticality is that it is well adapted to the probabilistic setting, since the functions involved are not rescaling of Schwartz functions, but of Fourier-randomized Schwartz functions.
\subsubsection{The parabolic case}\label{secpara} Start with (\ref{she}). Since we are interested in local theory, for simplicity we will periodize the time, so the noise
\[\xi\sim \sum_{(k,\ell)\in\mathbb{Z}^{d+1}}\langle k\rangle^{-\alpha}e^{i(k\cdot x+\ell\cdot t)}g_{k,\ell},\] where $g_{k,\ell}$ are i.i.d. Gaussians, and hence
\[u_{\mathrm{lin}}\sim \sum_{(k,\ell)\in\mathbb{Z}^{d+1}}\langle k\rangle^{-\alpha}(1+|k|^2+|\ell|)^{-1}e^{i(k\cdot x+\ell\cdot t)}g_{k,\ell}.\] Plugging into (\ref{she}) and assuming (say) $\mathcal{N}_p(u) = u^p$, we obtain that the second iteration
\[u_{\mathrm{iter}}\sim\sum_{(k,\ell)\in\mathbb{Z}^{d+1}}X_{k,\ell}e^{i(k\cdot x+\ell\cdot t)},\] where the random variables 
\[X_{k,\ell}:=(1+|k|^2+|\ell|)^{-1}\sum_{\substack{k_1+\cdots +k_p=k\\\ell_1+\cdots +\ell_p=\ell}}\prod_{q=1}^p\langle k_q\rangle^{-\alpha}(1+|k_q|^2+|\ell_q|)^{-1}g_{k_q,\ell_q}.\] Here, assuming all the $(k_q,\ell_q)$ involved are different, one may invoke the square root cancellation\footnote{This square root cancellation is derived from a classic large deviation property for the sum of independent Gaussian variables, as referenced in Lemma 4.4 in \cite {DNY20}.} for sums of independent Gaussians to obtain that with high probability,
\begin{equation}\label{sqrt1}|X_{k,\ell}|\lesssim(1+|k|^2+|\ell|)^{-1}\bigg(\sum_{\substack{k_1+\cdots +k_p=k\\\ell_1+\cdots +\ell_p=\ell}}\prod_{q=1}^p\langle k_q\rangle^{-2\alpha}(1+|k_q|^2+|\ell_q|)^{-2}\bigg)^{1/2}.\end{equation}
Now, assume we are dealing with \emph{high-high-to-high} interactions, then each $|k_q|\sim |k|\sim N$, and $|\ell_q|\sim|\ell|\sim N^2$,
hence by the heuristic dimensional counting argument of $(k_q,\ell_q)\in \mathbb{Z}^{d+1}$ for $q=1,...,p$ in the above sum,  the parabolic scaling being $d+2$ when $(k_q,\ell_q)\in \mathbb{Z}^{d+1}$,  the inequality \eqref{sqrt1} implies \[|X_{k,\ell}|\lesssim N^{-p\alpha-2p-2+(d+2)(p-1)/2}.\] From this it is easy to see that
\[\|u_{\mathrm{iter}}\|_{C_t^0C_x^s}\lesssim N^{-(p-1)\alpha-2+(d-2)(p-1)/2},\] where $s=\alpha+1-d/2$ is the regularity level of the solution. We thus see that the equation is
\[\left\{\begin{aligned}
&\mathrm{subcritical},&&\mathrm{if}\quad\alpha>\frac{d}{2}-1-\frac{2}{p-1}\quad\mathrm{or}\quad s>\frac{-2}{p-1},\\
&\mathrm{critical},&&\mathrm{if}\quad\alpha=\frac{d}{2}-1-\frac{2}{p-1}\quad\mathrm{or}\quad s=\frac{-2}{p-1},\\
&\mathrm{supercritical}, &&\mathrm{if}\quad\alpha<\frac{d}{2}-1-\frac{2}{p-1}\quad\mathrm{or}\quad s<\frac{-2}{p-1}.
\end{aligned}\right.\] In this way we get the critical threshold of (\ref{she}) is $C^s$ where $s=-2/(p-1)$, which agrees with the standard parabolic scaling as described in Hairer's work \cite{Hairer}.
\subsubsection{The hyperbolic case}\label{sec:hyperbolic-hhh} In the hyperbolic case (\ref{nlw}) we argue similarly, but also need to take into account the time oscillation. For simplicity we shall replace $u_1=\partial_tu(0)$ by $0$, which does not affect the scaling argument; then the linear component
\[u_{\mathrm{lin}}=\sum_{k\in\mathbb{Z}^d} \cos (t|k|) \langle k\rangle^{-\alpha-1} e^{ik\cdot x} g_k.\]
By using $\sin z = (e^{iz}-e^{-iz}/2)$, plugging the linear component into (\ref{nlw}) and assuming (say) $\mathcal{N}_p(u) = u^p$, we obtain that the second iteration
\[u_{\mathrm{iter}}\sim\sum_{k\in\mathbb{Z}^{d}}X_{k}e^{i(k\cdot x\pm_0 |k| t)},\] where the random variables 
\[X_{k}:=\langle k\rangle^{-1}\sum_{\substack{k_1+\cdots +k_p=k}}\frac{e^{it\Omega -1}}{\Omega}\prod_{q=1}^p \langle k_q \rangle^{-\alpha-1}g_{k_q},\] where \[\Omega =\mp_0 | k|\pm_1 |k_1|+\cdots \pm_p |k_p|.\] 
Here, assuming all the $k_q$ involved are different, one may invoke the square root cancellation for sums of independent Gaussians to obtain that with high probability,
\begin{equation}\label{sqrt2}|X_{k}|\lesssim \langle k\rangle^{-1}\bigg(\sum_{\substack{k_1+\cdots +k_p=k}}\langle \Omega\rangle^{-1}\prod_{q=1}^p\langle k_q\rangle^{-2\alpha-2}\bigg)^{1/2}.\end{equation}
Now, assuming we are dealing with \emph{high-high-to-high} interactions, then each $|k_q|\sim |k|\sim N$, hence by by omitting a loss of $N^\varepsilon$ and the lattice counting  with fixed $\Omega$ (Lemma \ref{counting:wave}), except for the case when $d\in \{1,2\}$ and $p=2$, we have 
\[|X_{k}|\lesssim N^{-p(\alpha+1)+d(p-1)/2-3/2}.\]  
From this it is easy to see that
\[\|u_{\mathrm{iter}}\|_{C_t^0H_x^s}\lesssim N^{-(p-1)\alpha+(d-2)(p-1)/2-3/2},\] where $s=\alpha+1-d/2$ is the regularity level of the solution. Also we hold the same bound in $C_t^0C_x^s$ by Khintchine's inequality (see the similar estimates in Section 7.2 of \cite{BDNY22}). We thus see that the equation is
\[\left\{\begin{aligned}
&\mathrm{subcritical},&&\mathrm{if}\quad\alpha>\frac{d}{2}-1-\frac{3}{2(p-1)}\quad\mathrm{or}\quad s>\frac{-3}{2(p-1)},\\
&\mathrm{critical},&&\mathrm{if}\quad\alpha=\frac{d}{2}-1-\frac{3}{2(p-1)}\quad\mathrm{or}\quad s=\frac{-3}{2(p-1)},\\
&\mathrm{supercritical}, &&\mathrm{if}\quad\alpha<\frac{d}{2}-1-\frac{3}{2(p-1)}\quad\mathrm{or}\quad s<\frac{-3}{2(p-1)}.
\end{aligned}\right.\] In this way we get the critical threshold of (\ref{nlw}) is $C^s$ where $s=-\frac{3}{2(p-1)}$ except for the case when $d=1, 2$ and $p=2$.

In the case $d=1$ and $p=2$, by Lemma \ref{counting:wave}  we have $|X_{k}|\lesssim N^{-p(\alpha+1)+d(p-1)/2-1}.$ And hence it is easy to see that the critical threshold is $C^s$ where $s=-1$.
In the case $d=2$ and $p=2$, by Lemma \ref{counting:wave}  we have $|X_{k}|\lesssim N^{-p(\alpha+1)+d(p-1)/2-1}.$ And hence it is easy to see that the critical threshold is $C^s$ where $s=-\frac{5}{4}$.

\subsubsection{The dispersive case}\label{secdisper}
In the dispersive case (\ref{nls}), then the linear component
\[u_{\mathrm{lin}}=\sum_{k\in\mathbb{Z}^d}  \langle k\rangle^{-\alpha-1} e^{i(k\cdot x+|k|^2t)} g_k.\]
Plugging into (\ref{nls}) and assuming (say) $\mathcal{N}_p(u)= u^{\pm_1}\cdots u^{\pm_p}$  (denote that $u^+:=u$ and $u^-:=\overline{u}$), we obtain that the second iteration
\[u_{\mathrm{iter}}\sim\sum_{k\in\mathbb{Z}^{d}}X_{k}e^{i(k\cdot x\pm_0 |k|^2 t)},\] where the random variables 
\[X_{k}:=\sum_{\substack{\pm_1 k_1+\cdots \pm_p k_p=k}}\frac{e^{it\Omega -1}}{\Omega}\prod_{q=1}^p \langle k_q \rangle^{-\alpha-1}g_{k_q},\] where \[\Omega =\mp_0 | k|^2\pm_1 |k_1|^2+\cdots \pm_p |k_p|^2.\] 
Here, assuming there are no \emph{parings}\footnote{$k_i
\neq k_j$ if the corresponding  signs $\pm_i$ and $\pm_j$ are the opposite.} in the $p$-th order power $\mathcal{N}_p$, one may invoke the square root cancellation for sums of independent Gaussians to obtain that with high probability,
\begin{equation}\label{sqrt3}|X_{k}|\lesssim \bigg(\sum_{\substack{\pm_1 k_1+\cdots \pm_p k_p=k}}\langle \Omega\rangle^{-1}\prod_{q=1}^p\langle k_q\rangle^{-2\alpha-2}\bigg)^{1/2}.\end{equation}
Now, assuming we are dealing with \emph{high-high-to-high} interactions, then each $|k_q|\sim |k|\sim N$. 
 From this it is easy to see that
\begin{multline}\label{est:nls}
\|u_{\mathrm{iter}}\|_{C_t^0H_x^s}\lesssim N^{-(p-1)(\alpha+1)-d/2} \times\\  \left(\# \{k,k_q\in \mathbb Z^d \text{ for all }q: |k|\sim N, |k_q|\sim N, \text{ for all }q, \pm_1 k_1+\cdots \pm_p k_p=k, \Omega\text{ is fixed}\}\right)^{\frac{1}{2}}\end{multline}
 where $s=\alpha+1-d/2$ is the regularity level of the solution. 
Hence by omitting a loss of $N^\varepsilon$ and the lattice counting  with fixed $\Omega$ (Lemma \ref{counting:nls}), we have
\[\|u_{\mathrm{iter}}\|_{C_t^0H_x^s}\lesssim N^{-(p-1)\alpha+(d-2)(p-1)/2-1}.\]
Also we hold the same bound in $C_t^0C_x^s$ by Khintchine's inequality.
We thus see that the equation is
\[\left\{\begin{aligned}
&\mathrm{subcritical},&&\mathrm{if}\quad\alpha>\frac{d}{2}-1-\frac{1}{p-1}\quad\mathrm{or}\quad s>\frac{-1}{p-1},\\
&\mathrm{critical},&&\mathrm{if}\quad\alpha=\frac{d}{2}-1-\frac{1}{p-1}\quad\mathrm{or}\quad s=\frac{-1}{p-1},\\
&\mathrm{supercritical}, &&\mathrm{if}\quad\alpha<\frac{d}{2}-1-\frac{1}{p-1}\quad\mathrm{or}\quad s<\frac{-1}{p-1}.
\end{aligned}\right.\] In this way we get the critical threshold of (\ref{nls}) is $C^s$ where $s=-\frac{1}{p-1}$.

\begin{rem}
The estimate in \eqref{est:nls} heavily relies on the square root cancellation for sums of independent Gaussians. If we replace the Gaussians in the data \eqref{data3} simply by $1$ and restrict the data to frequency $\sim N$, then the data would be a smooth data whose Fourier modes are supported around $N$. With such smooth deterministic data, we lose the square root gain in \eqref{est:nls}, and by the similar calculations as above we get another critical threshold $s= \frac{d}{2}-\frac{2}{p-1}$ of \eqref{nlw} in $H^s$, which coincides the usual deterministic scaling for \eqref{nls}.
Similarly in the hyperbolic case with such smooth deterministic data, we can also obtain the usual deterministic scaling critical $s$ for \eqref{nlw}.
\end{rem}

\section{Discussions}\label{estimates}
\subsection{Possible discrepancies}\label{sec:discrep} We should point out that the scaling heuristics in Section \ref{scaling} provide only a \emph{guiding principle} to the probabilistic well-posedness problem of the corresponding dynamics and \emph{should not be understood} that the actual threshold between a.s. local well-posedness and ill-posedness is always given by scaling criticality. In most cases this is indeed true, including for example general NLS equations of cubic or higher nonlinearity \cite{DNY20}, but in some cases, especially concerning low dimensions and/or low degree nonlinearity, discrepancies may occur between the scaling prediction and well-posedness result\footnote{For example \cite{EXTRA1,EXTRA2} focus on such discrepancies.}.

Such discrepency is not uncommon in other settings involving the notion of scaling, in fact it is well-known that:
\begin{itemize}
\item In the deterministic setting, the local well-posedness threshold for NLW equations generally does not equal that predicted by scaling criticality due to the Lorentzian symmetry; similar phenomena happen for NLS equations involving negative regularity data.
\item In the case of stochastic heat equations, the local well-posedness threshold again deviates from the parabolic scaling prediction, if the noise involved is rougher than the spacetime white noise (for a more comprehensive description see  \cite{MoWe} and references
therein).
\end{itemize}

In our case, the discrepancy is mainly caused by two mechanisms: (1) the high-high-to-low interaction, which is the same reason for the discrepancy in the stochastic heat equations, and (2) the anomalies occurring in various counting estimates depending on the specific dispersion relation, which is related to the discrepancy in deterministic problems for NLW. We now describe these mechanisms in more detail.

\subsubsection{The role of high-high-to-low interactions}\label{sec:hhl}
As indicated in Section \ref{sec:criticality}, the probabilistic scaling criticality is determined through heuristic calculations of the \emph{high-high-to-high} interactions, which are also taken into account when obtaining the critical exponent for the usual deterministic scaling. As indicated in Section \ref{sec:discrep}, the \emph{high-high-to-low} interactions may also play a role in the local theory of the random data problem for the dispersive and wave equations in some special cases which we will discuss below.

Following similar calculations as in Section \ref{sec:criticality} (see \eqref{sqrt1} for the heat equation, \eqref{sqrt2} for the wave equation and \eqref{sqrt3} for the Schr\"odinger equation) and assuming we are dealing with \emph{high-high-to-low} interactions (each $|k_q|\sim N$ and $|k|\sim 1$) then up to an $N^\varepsilon$ loss and using the lattice counting lemmas in Section \ref{sec:counting}, we can obtain that
\[\left\{\begin{aligned}
&\text{heat:} &&\|u_{\mathrm{iter}}\|_{C_t^0C_x^s}\lesssim N^{-sp-(d+2)/2},\\
&\text{wave:}&&
\|u_{\mathrm{iter}}\|_{C_t^0C_x^s}\lesssim\begin{cases} N^{-sp-d/2} \quad &\text{when }\mathcal{N}_p(u)=u^2,\\N^{-sp-(d+1)/2} \quad &\text{otherwise,}\end{cases}\\
&\text{Schr\"odinger:} &&\|u_{\mathrm{iter}}\|_{C_t^0C_x^s}\lesssim \begin{cases} N^{-sp-(d+1)/2}\quad &\text{when }\mathcal{N}_p(u)=|u|^2,\\N^{-sp-(d+2)/2}\quad &\text{otherwise}.\end{cases},
\end{aligned}\right.\] where $s=\alpha+1- d/2$.
As a result, we can see the critical threshold arising from the high-high-to-low interactions is $C^{s}$ where, 
\[\left\{\begin{aligned}
&\text{heat:} &&s=-\frac{d+2}{2p},\\
&\text{wave:}&&
s=\begin{cases} -\frac{d}{2p}  \quad &\text{when }\mathcal{N}_p(u)=u^2,\\-\frac{d+1}{2p}  \quad &\text{otherwise,}\end{cases}\\
&\text{Schr\"odinger:} && s=\begin{cases}-\frac{d+1}{2p}\quad &\text{when }\mathcal{N}_p(u)=|u|^2,\\-\frac{d+2}{2p}\quad &\text{otherwise}.\end{cases}
\end{aligned}\right.\]

Let us compare the \emph{probabilistic scaling} critical threshold $s_{pr}$ and the critical threshold $s$ above which we denote now by $s_{hhl}$ for the \emph{high-high-to-low interactions} in the context of the stochastic heat equation \eqref{she}, the semilinear wave equation \eqref{nlw} and the semilinear Schr\"odinger equation \eqref{nls}. As shown in Table \ref{intro:scalingfig3}, for \eqref{she} we obtain $s_{pr}\geq  s_{hhl}$ when $d>2$ and $p\geq \frac{d+2}{d-2}$; for \eqref{nlw} we obtain $s_{pr}\geq  s_{hhl}$ when 
$d>2$, $p\geq \max (\frac{d+1}{d-2}, 3)$  and when $d\geq 5$, $p=2$; for \eqref{nls}  we obtain $s_{pr}\geq  s_{hhl}$   when $d\geq 3$, $\mathcal{N}_p(u)=|u|^2$ and when $p\geq \frac{d+2}{d}$, $\mathcal{N}_p(u)\neq |u|^2$.
\renewcommand{\thefootnote}{\fnsymbol{footnote}}

\begin{table}\renewcommand{\arraystretch}{2.3}
\begin{tabular}{ |p{2.5cm}<{\centering} |p{1.5cm}<{\centering}|p{2cm}<{\centering}|p{2cm}<{\centering} | } 
 \hline
  & $\mathcal{N}_p(u)$& $s_{pr}$& $s_{hhl}$ \\ \hline
Heat  & all $p$  & $-\frac{2}{p-1}$ &$-\frac{d+2}{2p}$   \\ \hline
 Wave & $p=2$& $-\frac{3}{2}$\tablefootnote[1]{{Note that, for $d\in\{1,2\}$, the actual threshold is higher than this value $-3/2$, which has to do with specific counting estimates for wave equations; see Lemma \ref{counting:wave}. This same discrepancy between scaling critical and well-posedness result also happens for the deterministic wave equation in 2D.}}  &$-\frac{d}{4}$  \\ 
 & $p\geq 3$& $-\frac{3}{2(p-1)}$  &$-\frac{d+1}{2p}$  \\ 
 \hline
 Schr\"odinger & $p=2$\tablefootnote[2]{Here $p=2$ only means the case $\mathcal{N}_p(u)=|u|^2$. When $\mathcal{N}_p(u)=u^2$ or $\overline{u}^2$, $s_{hhl}$ still follows $-\frac{d+2}{2p}=-\frac{d+2}{4}$}  & $-1$ &$-\frac{d+1}{4}$ \\ 
 & $p\geq 3$& $-\frac{1}{p-1}$  &$-\frac{d+2}{2p}$  \\ 
 \hline
\end{tabular}
\caption[hhh]{\small Probabilistic scaling criticality v.s. high-high-to-low criticality}
\label{intro:scalingfig3}
\end{table}

\renewcommand{\thefootnote}{\arabic{footnote}}

\subsubsection{The role of dispersion relations  and counting estimates}\label{sec:counting}
The scaling estimates in Sections \ref{sec:criticality} and \ref{sec:hhl} crucially rely on the following counting estimates which are intimately related to the dispersion at hand.
Although in Lemma \ref{counting:wave}
and Lemma \ref{counting:nls} we only state upper bounds, in most cases these are optimal. For the sake of brevity, we do not repeat the proofs in this note, but the reader may find these in \cite{BDNY22, Dean,  DNY19, DNY20}. {Note in particular the less favorable counting estimate in (\ref{counting:eq1}) when $A\ll N$, which corresponds to high-high-to-low interactions discussed above, and the different counting estimates (\ref{counting:eq2}) and (\ref{counting:eq3}), which has to do with the anomalous scaling for wave equation in 2D.}

\begin{lem}[The counting lemma for wave equations]\label{counting:wave}
Given  dyadic numbers $A, N$ and $a\in \mathbb{Z}^d$ satisfying $|a|\sim A$, we have the following lattice point counting bounds
\begin{enumerate} 
\item  when $d\geq 3$ \begin{equation}\label{counting:eq1}
\sup_{m\in \mathbb{Z}} \# \{n\in \mathbb{Z}^d: |n|\sim N, |\langle a+n\rangle-\langle n\rangle -m|\leq 1\} \lesssim \min (A, N)^{-1}  N^d.
\end{equation}
\item   when $d\geq 3$ \begin{equation}\label{counting:eq2}
\sup_{m\in \mathbb{Z}} \# \{n\in \mathbb{Z}^d: |n|\sim N, |\langle a+n\rangle+\langle n\rangle -m|\leq 1\} \lesssim N^{d-1}.
\end{equation}
\item  when $d=2$
\begin{equation}\label{counting:eq3}
\sup_{m\in \mathbb{Z}} \# \{n\in \mathbb{Z}^2: |n|\sim N, |\langle a+n\rangle -\langle n\rangle -m|\leq 1\} \lesssim N^{\frac{3}{2}}.
\end{equation}
\item  when $d=1$
\begin{equation}\label{counting:eq3}
\sup_{m\in \mathbb{Z}} \# \{n\in \mathbb{Z}: |n|\sim N, |\langle a+n\rangle -\langle n\rangle -m|\leq 1\} \lesssim N.
\end{equation}
\end{enumerate}
\end{lem}

\begin{lem}[The counting lemma for Schr\"odinger equations]\label{counting:nls}
Given  dyadic numbers $N$ and $a\in \mathbb{Z}^d$, we have the following lattice point counting bounds
\begin{enumerate} 
\item \begin{equation}\label{counting:eq4}
\sup_{m\in \mathbb{Z}} \# \{n\in \mathbb{Z}^d: |n|\sim N, \big| |a+n|^2-| n|^2 -m\big|\leq 1\} \lesssim  N^{d-1}.
\end{equation}
\item \begin{equation}\label{counting:eq4}
\sup_{m\in \mathbb{Z}} \# \{n\in \mathbb{Z}^d: |n|\sim N, \big| |a+n|^2+| n|^2 -m\big|\leq 1\} \lesssim  N^{d-2+\varepsilon}.
\end{equation}
\end{enumerate}
\end{lem}

\subsection{The comparison with existing results and with Gibbs measures}\label{gibbs}
The Gibbs measure associated to the stochastic heat equation with space-time white noise \eqref{she}, the nonlinear wave equation \eqref{nlw} and the nonlinear Schr\"odinger equation \eqref{nls} 
is known to exist for $d=1,2$ and for all power nonlinearity $\mathcal{N}_p(u)$ with odd $p\geq 3$ (see \cite{Bourgain94, Nelson66} for details).  When $d=3$, the Gibbs measure is also known to exist but only when $p=3$; this is so-called $\Phi^4_3$ measure (see \cite{BaGu0, BaGu, GJ1} and references therein for details). However, one has the marginal triviality of the $\Phi^4_d$-measure for $d\geq 4$ (see the works of Aizenman \cite{A81}, Fr\"{o}hlich \cite{F82}, and Aizenman-Duminil-Copin \cite{AC21}).  Loosely speaking, this implies that the $\Phi^4_d$-measure in dimension $d\geq 4$  essentially yields a Gaussian measure, for any renormalization of the potential energy term  giving rise to a well defined measure.

To rigorously justify the invariance of the Gibbs measures under the corresponding dynamics, we need to first fully understand the a.s. local wellposedness of \eqref{she}, \eqref{nlw} and \eqref{nls} with Gaussian data ($\alpha=0$) \eqref{data1}, \eqref{data2} and \eqref{data3} respectively. The Gibbs measure is supported in $H^{s_G-}$ where $s_G = 1-\frac{d}{2}$.

An interesting question arises in the context of dispersive and wave equations. Let us focus on the NLS. In 1996, Bourgain proved in a seminal paper \cite{Bourgain} the a.s. global wellposedness and the invariance of the Gibbs measure under the dynamics of the 2D cubic NLS. But since the measure is known to exist in 2D for odd $p\geq 5$ and in 3D for $p=3$, it is natural  to ask about the a.s. global wellposedness and the invariance of the Gibbs measure in these cases.
See also \cite{Bourgain94, LRS} for more works about the invariance of Gibbs measures under nonlinear Sch\"odinger equations.
 In our paper \cite{DNY19}, we indeed proved these facts for $d=2$ and odd $p\geq 5$. 
The interesting question is why is it that the random averaging operator method introduced in \cite{DNY19} does not suffice when $d=3$ and $p=3$.
Note that the deterministic scaling $s_{cr} = \frac{d}{2} -\frac{2}{p-1}$ when $d=2$ and $p=5$ is $H^{\frac{1}{2}}(\mathbb T^d)$ same as when $d=3$ and $p=3$.
What tells these two cases apart is precisely the probabilistic scaling $s_{pr}$.
When  $d=2$ and $p=5$ the probabilistic scaling $s_{pr}=-\frac{1}{4}< 0- = s_{G}-$, so the the a.s. local wellposedness problem corresponding to the Gibbs measure is probabilistic subcritical. But when $d=3$ and $p=3$ 
 the probabilistic scaling $s_{pr}=-\frac{1}{2}=  s_{G}$, so the the corresponding a.s. local wellposedness problem is probabilistic critical.
The probabilistic scaling is thus the correct lens through which we should view the probabilistic wellposedness problems.

In the case of the cubic nonlinear wave equation in dimension three (the hyperbolic $\Phi^4_3$ model), the statistiscal ensemble of its  associated Gibbs measure is probabilistical subcritical. In \cite{BDNY22} in joint work with B. Bringmann we prove the invariance of the Gibbs measure under the dynamics of the three-dimensional cubic wave equation (see \cite{B99, F85, GKO18, OT20, OOT21, Z94} for more about the hyperbolic $\Phi^{p+1}$ models). This result is the hyperbolic counterpart to seminal works on the parabolic $\Phi^4_3$ model by Hairer \cite{Hairer} and Hairer-Matetski \cite{HM18}.  See \cite{AK20,CC18, DPD03, GH21,GIP15, I87, MoWe, MouWe, MouWe2, PW81} and references therein for more works on the parabolic $\Phi^{p+1}$ models.

To end this section,  Table \ref{intro:scalingfig1} shows the optimal number of derivatives $\beta$ that the initial data can be rougher than the Gaussian data ($\alpha=0$) corresponding to the Gibbs measure, when analyzing the high-high into high (probabilistical scaling) interactions. In contrast in Table \ref{intro:scalingfig2} we lay out the optimal number of
 derivatives $\beta$ that arises when analyzing the high-high into low interactions.
\renewcommand{\thefootnote}{\fnsymbol{footnote}}

\begin{table}\renewcommand{\arraystretch}{2.0}
\begin{tabular}{ |p{2cm}<{\centering}|p{2cm}<{\centering}|p{2cm}<{\centering}|p{2cm}<{\centering}|p{2cm}<{\centering}|p{2cm}<{\centering}|p{2cm}<{\centering} | } 
 \hline
  & 2D Quad.& 2D Cubic &3D Quad. & 3D Cubic& 4D Quad.& 4D Cubic\\ \hline
Heat  & \cellcolor{green!10} $\beta= 2$ & \cellcolor{green!10}$\beta=1$ &\cellcolor{green!10}$\beta=\frac{3}{2}$   & \cellcolor{green!10}$\beta=\frac{1}{2}$ & \cellcolor{green!10}$\beta=1$ & \cellcolor{red!10}$\beta=0$ \\ \hline
 Wave & \cellcolor{green!10}$\beta=\frac{5}{4}$&\cellcolor{green!10} $\beta=\frac{3}{4}$  &\cellcolor{green!10}$\beta=1$  &\cellcolor{green!10} $\beta=\frac{1}{4}$ & \cellcolor{green!10}$\beta=\frac{1}{2}$ & \cellcolor{red!20}$\beta=-\frac{1}{4}$ \\ 
 \hline
Schr\"{o}dinger  &\cellcolor{green!10} $\beta=1$ &\cellcolor{green!10} $\beta=\frac{1}{2}$  &\cellcolor{green!10}$\beta=\frac{1}{2}$  &\cellcolor{red!10} $\beta=0$ &\cellcolor{red!10} $\beta=0$ & \cellcolor{red!20} $\beta=-\frac{1}{2}$  \\ 
 \hline
\end{tabular}
\caption[hhh]{\small $\beta$ is the number of the derivatives that the initial data can be rougher than the Gaussian data $\alpha=0$ corresponding to the Gibbs measure concerning the high-high into high (probabilistical scaling) interaction.}
\label{intro:scalingfig1}
\end{table}

\renewcommand{\thefootnote}{\arabic{footnote}}

\begin{table}\renewcommand{\arraystretch}{2.0}
\begin{center}\begin{tabular}{ |p{2cm}<{\centering}|p{2cm}<{\centering}|p{2cm}<{\centering}|p{2cm}<{\centering}|p{2cm}<{\centering}|p{2cm}<{\centering}|p{2cm}<{\centering} | } 
 \hline
  & 2D Quad.& 2D Cubic &3D Quad. & 3D Cubic& 4D Quad.& 4D Cubic\\ \hline
Heat &\cellcolor{green!10} $\beta= 1$ & \cellcolor{green!10}$\beta=\frac{2}{3}$ &\cellcolor{green!10}$\beta=\frac{3}{4}$   & \cellcolor{green!10}$\beta=\frac{1}{3}$ & \cellcolor{green!10}$\beta=\frac{1}{2}$ & \cellcolor{red!10}$\beta=0$ \\ \hline
Wave & \cellcolor{green!10}$\beta=\frac{1}{2}$ &\cellcolor{green!10} $\beta=\frac{1}{2}$  &\cellcolor{green!10}$\beta=\frac{1}{4}$  &\cellcolor{green!10} $\beta=\frac{1}{6}$ & \cellcolor{red!10}$\beta=0$ &\cellcolor{red!20} $\beta=-\frac{1}{6}$ \\ 
 \hline
Schr\"{o}dinger  &  \cellcolor{green!10}$\beta=\frac{3}{4}$ & \cellcolor{green!10} $\beta=\frac{2}{3}$  & \cellcolor{green!10}$\beta=\frac{1}{2}$  &  \cellcolor{green!10}$\beta=\frac{1}{3}$ &  \cellcolor{green!10}$\beta=\frac{1}{4}$ &  \cellcolor{red!10}$\beta=0$  \\ 
 \hline
\end{tabular}
\end{center}
\caption{\small $\beta$ is the number of the derivatives that the initial data can be rougher than the Gaussian data ($\alpha=0$) corresponding to the Gibbs measure concerning the high-high into low interaction.}
\label{intro:scalingfig2}
\end{table}


\subsection{Long time solutions and connection to wave turbulence theory}\label{WTT} The same scaling heuristics in Section \ref{scaling} can also be applied to more general context beyond local well-posedness. For example, fix a dyadic scale $N$, and assume that the initial data is supported at scale $N$ as in Section \ref{sec:criticality}. Suppose instead of local well-posedness (i.e. existence up to time $1$) we are interested in existence up to time scale $T\sim N^\beta$ for some $\beta$, where $\beta>0$ corresponds to long-time existence result as $N\to\infty$, and $\beta<0$ corresponds to short-time existence result as $N\to\infty$. Then, we may consider the following questions:
\begin{itemize}
\item Suppose $\beta$ is fixed, what is the optimal decay rate of the initial data, in terms of the exponent $\alpha$ in Sections \ref{secpara}--\ref{secdisper}, such that solution to the given equation exists with high probability up to time $T=N^\beta$?
\item Conversely, given the decay rate $\alpha$, what is the best time $T=N^\beta$ up to which a solution exists with high probability?
\end{itemize}

In fact, the above questions are exactly those encountered in the mathematical theory of \emph{wave turbulence}.
The key prediction of the wave turbulence theory, starting from random initial data (or stochastically forced) problems, is the so-called \emph{wave kinetic equation}, which not only claims the existence of solutions, but also calculates the effective dynamics of statistical quantities of the solution. This effective dynamics occurs at a particular time scale $T_{\mathrm{kin}}$, called the \emph{kinetic} or \emph{Van Hove} time scale.

In fact, this notion of $T_{\mathrm{kin}}$ is precisely the time scale predicted by the heuristic calculations in Section \ref{secpara}--\ref{secdisper}: given the decay of initial data quantified by the exponent $\alpha$, this $T_{\mathrm{kin}}=N^\beta$ is exactly the value such that the ($C^s$ or) $H^s$ norms of $u_{\mathrm{iter}}$ and $u_{\mathrm{lin}}$ are comparable at time $T$. In particular, the probabilistically critical, subcritical and supercritical problems in the sense of Section \ref{scaling}, precisely correspond\footnote{In wave turbulence theory it is customary to perform another reduction so the space scale becomes $N$ and time scale becomes $N^2$; in this setting the probabilistically critical problem would correspond to $T_{\mathrm{kin}}=N^2$.} to those problems where $T_{\mathrm{kin}}=1$, $T_{\mathrm{kin}}\gg1$ and $T_{\mathrm{kin}}\ll1$.

In the paper \cite{DNY20}, as a consequence of 
the sharp subcritical a.s. local wellposedness of \eqref{nls} in $H^s(\mathbb T^d)$ with $s>s_{pr}$,
we also obtained the long time existence\footnote{Thus, in particular, proving that with high probability, there is no energy cascade between Fourier modes (i.e. $|\widehat{u}(t,k)|^2\approx|\widehat{u}(0,k)|^2$ with negligible error for large $N$).} of solutions for  well prepared smooth random data  --such as that arising in the derivation of the wave kinetic equation as we explained above  up to the time $T=N^{(p-1)(s-s_{pr})-}.$
In the case of general odd $p$, one can show that the kinetic time scale $T_{kin}$ refered above is $N^{2(p-1)(s-s_{pr})+2}$.
After a suitable rescaling, the time $T$ obtained in \cite{DNY20} reaches $N^{-\varepsilon} T_{kin}$ when $s-s_{pr}=\varepsilon/(p-1)$, for all $p\geq 3$
odd.

 The physical theory of wave turbulence started in the 1920s \cite{Pei} and had developed into a substantial field of research with a wide range of scientific and practical applications; the rigorous mathematical theory started much later around the 2000s, and has become particularly active in recent years, see \cite{CG, CG20, DH,DH2} and references therein.
It has been a major open problem in wave turbulence theory to derive the wave kinetic equation up to the time scale $T_{\mathrm{kin}}$; recently this has been accomplished by the first author with Z. Hani in the work \cite{DH2} for the cubic NLS, and the methods used in \cite{DH2} are closely related to those employed in \cite{DNY19,DNY20}. Note that \cite{DH2} covers the probabilistically critical case $T_{\mathrm{kin}}=1$, so it in fact solves a \emph{critical} problem, which may shed some light on the Gibbs measure invariance problem described in Section \ref{gibbs} above. Though, this latter problem is still much harder, due to the many potential logarithmic divergences.

\subsection{Other geometries} 
In the context of expanding to more general geometries, including $\mathbb{T}^d$, we encounter different scenarios.
In the case of compact manifolds, 
the canonical randomization would be based on the spectral expansion of the Laplacian. Here, the probabilistic scaling depends on the \emph{global geometry} of the underlying manifold. This is because of the fact that the randomized data is \emph{not} localized and has the same amplitude at each point within the domain. In comparison, the \emph{deterministic scaling} threshold remains independent of the geometry, because it corresponds to the data zoomed out at a point, which is localized. On the other hand, although the gaussian measure can be defined canonically on manifolds by expanding the eigenfunctions of the Laplacian and this concept itself does not depend too much on the geometry, when we treat nonlinear problems, the (probabilistic) analysis in this case will depend on multilinear integrals involving eigenfunctions which depend on the underlying geometry.

In the case of non-compact manifolds, for instance $\mathbb{R}^d$, the canonical randomization based on Laplacian eigenfunctions (i.e., $e^{i\xi\cdot x}$) leads to initial data with infinite $L^2$ mass, which is not feasible\footnote{However, it is compatible with the wave equation due to the finite speed of propagation; specifically, the result of \cite{BDNY22} is expected to hold also for the Gibbs measure on $\mathbb{R}^3$. We refrain from discussing this here, but for more details, refer to \cite{BDNY22}.} for the Schr\"odinger equation (\ref{nls}). There is an alternative ``Wiener randomization" approach that involves sectioning the Fourier space into unit boxes and performing randomization within each box (refer to, for example, \cite{ZhF, LM, BOP}). This randomization generates localized initial data, which is not conserved by the linear flow. In the $\mathbb{R}^d$ scenario, this results in a critical threshold $s_{p}=-3/(4r)$ that is lower than that of $\mathbb{T}^d$.

  \end{document}